\documentclass[12pt,leqno,amstex]{article}
\usepackage{amsfonts}
\usepackage{amsfonts,amssymb,amsmath}
\newcommand{\ve}{\varepsilon }
\newcommand{\R}{\mathbb{R}}

\begin{document}
\title{Spatially Inhomogeneous Bernstein's  Problem and De Giorgi's Conjecture}
\author{By {\sc Bendong Lou}
}
\date{}
\maketitle
\begin{center}

\begin{minipage}{120mm}
{{\bf Abstract}.\ \ In this note, we propose Bernstein's  problem
and De Giorgi's conjecture for spatially inhomogeneous equations, as
well as De Giorgi's conjecture for system of reaction-diffusion
equations.}
\end{minipage}
\end{center}
\ \\

 \baselineskip 18pt
\section{Bernstein's  problem and De Giorgi's conjecture}

 Consider
\begin{equation}\label{RD-homo}
u_t = \Delta u + \frac{1}{\ve^2}(u - u^3) , \ \ \ \ \ \ \ {\rm
\mbox{\boldmath $x$}} = ({\rm \mbox{\boldmath $x$}}', x_{n}) \in
\mathbb{R}^{n},\ \ t>0,
\end{equation}
where $ {\rm \mbox{\boldmath $x$}}'  = (x_1, \cdots, x_{n-1})$. It
was shown in \cite{Chen, nmhs} (and reference therein) that, as
$0<\ve\ll 1$, a sharp internal layer (interface) of $u$ develops
between the two regions $\{u\approx 1\}$ and $\{u\approx -1\}$ in
time scale $\tau = t/\ve^2$. In time scale $t$, the interface
propagates according to the following mean curvature flow equation:
\begin{equation}\label{MCF-homo}
V= -(n-1)\kappa ,\ \ \ \ \ \ \ \ {\rm \mbox{\boldmath $x$}} \in
\varGamma_t, \ t>0,
\end{equation}
where $\varGamma_t := \{{\rm \mbox{\boldmath $x$}}\in \R^{n}\ |\ u
({\rm \mbox{\boldmath $x$}}, t) = 0\}$ is a hypersurface, $V$
denotes the normal velocity of $\varGamma_t$ and $\kappa$ denotes
the mean curvature.

The well known Bernstein's  problem is about the steady state of
\eqref{MCF-homo}:

\vskip 6pt \underline{\bf Bernstein's  Problem}: {\it Let
$h:\R^{n-1}\rightarrow \R$ be a $C^2$ function,
$\mathfrak{S}(h)\subset \R^{n}$ be the graph of $h$. If $\kappa =0$
on $\mathfrak{S}(h)$, then $\mathfrak{S}(h)$ is a hyperplane. }
\vskip 6pt

The answer was shown to be positive for $n\leq 8$ (E. De Giorgi,
Almgren, Simons) and negative for $n \geq 9$ (E. Bombieri -De Giorgi
- E. Giusti). Clearly, Berstein's problem is related to the shape of
the stationary solution of \eqref{RD-homo}. For which, De Giorgi
\cite{DeG} proposed the following conjecture:

 \vskip 6pt \underline{\bf De Giorgi's Conjecture}: {\it Assume that $u$ is
 a stationary, entire solution
(defined for all ${\rm \mbox{\boldmath $x$}} \in \R^{n}$) of
\eqref{RD-homo}, that it satisfies $|u| \leq 1$ and ${\frac
{\partial u } {\partial x_n}} >0$. Then, at least for $n\leq 8$, the
level sets of $u$ must be hyperplanes.}

\vskip 6pt This conjecture was partially proved recently in
\cite{GG}, \cite{Sa} and references therein.

\section{Spatially Inhomogeneous Bernstein's  Problem and De
Giorgi's Conjecture}

In the last two decades, many anthors studied the following
spatially inhomogeneous reaction-diffusion equation
\begin{equation}\label{RD-inhomo}
u_t = \nabla ( A ({\rm \mbox{\boldmath $x$}}) \nabla u ) + {\frac
{1}{\varepsilon^2}} B ({\rm \mbox{\boldmath $x$}})
 u \left( Z^2 ({\rm \mbox{\boldmath $x$}})  - u^2 \right),
 \ \ \ \ \ \ \ {\rm \mbox{\boldmath $x$}}  \in \mathbb{R}^{n},\ t>0,
\end{equation}
where $A, B$ and $Z$ are bounded, smooth functions with positive
infimums. A special case is $Z\equiv 1$ (typical double-well
potential with equal-well-depth).

It is easily seen that, when $0< \varepsilon \ll 1$,
\eqref{RD-inhomo} has exactly three stationary solutions: $Z_+ ,\ 0$
and $Z_-$ with $Z_\pm \thickapprox \pm Z$. As the homogeneous case,
it was shown in \cite{HKMN, Lou, nmhs} that, in time scale $t$, the
law of the propagation of the level set $\varGamma_t := \{{\rm
\mbox{\boldmath $x$}}\ |\ u ({\rm \mbox{\boldmath $x$}}, t) = 0\}$
is
\begin{equation}\label{Vg}
V= -(n-1) a( {\rm \mbox{\boldmath $x$}} ) \kappa + c( {\rm
\mbox{\boldmath $x$}} ) \nabla d( {\rm \mbox{\boldmath $x$}}) \cdot
{\rm \mbox{\boldmath $n$}}  \qquad \qquad {\rm for} \quad {\rm
\mbox{\boldmath $x$}} \in \varGamma_t
\end{equation}
where \mbox{\boldmath $n$} is the normal direction to $\varGamma_t$,
$V$ and $\kappa$ are as above, $a,c,d$ are bounded, smooth functions
with $\inf a >0$. Similar to the homogeneous case, we propose the
following problem:

\vskip 6pt \underline{\bf Spatially Inhomogeneous Bernstein's
Problem}: {\it Let $h:\R^{n-1}\rightarrow \R$ be a $C^2$ function,
$\mathfrak{S}(h)\subset \R^{n}$ be the graph of $h$. Assume that
\begin{equation}\label{MCF-inhomo-1}
-(n-1) a( {\rm \mbox{\boldmath $x$}} ) \kappa + c( {\rm
\mbox{\boldmath $x$}} ) \nabla d( {\rm \mbox{\boldmath $x$}}) \cdot
{\rm \mbox{\boldmath $n$}} =0 \qquad \qquad {\rm for} \quad {\rm
\mbox{\boldmath $x$}} \in \mathfrak{S}(h).
\end{equation}
Then there exists $h_0:\R^{n-1}\rightarrow \R$, whose graph is a
hyperplane, such that

{\rm (i)} $h-h_0$ is quasi-periodic if $a,c,d$ are periodic;

{\rm (ii)} $h- h_0$ is almost periodic if $a,c,d$ are almost
periodic.}

\vskip 6pt

In almost periodic case (ii), a primary analysis indicates that an
additional condition maybe also needed: For any ball $B \subset
\R^n$,
\begin{equation}\label{ap-int-1}
\int_{B} [a({\rm \mbox{\boldmath $x$}}) - a_0 ] d{\rm
\mbox{\boldmath $x$}}, \quad \int_{B} [c({\rm \mbox{\boldmath $x$}})
- c_0 ] d{\rm \mbox{\boldmath $x$}}, \quad \int_{B} [d({\rm
\mbox{\boldmath $x$}}) - d_0 ] d{\rm \mbox{\boldmath $x$}}
\end{equation}
are bounded by $M$ (independent of $B$), where $a_0, c_0$ and $d_0$
are the averages of $a,c$ and $d$, respectively.

\vskip 6pt

 \underline{\bf Spatially Inhomogeneous De Giorgi's Conjecture}.
{\it Assume that $u$ is an entire, stationary solution of
\eqref{RD-inhomo}, that it satisfies
\begin{equation}\label{2limits}
Z_- ({\rm \mbox{\boldmath $x$}}) \leq u ({\rm \mbox{\boldmath $x$}})
\leq Z_+ ({\rm \mbox{\boldmath $x$}}),\ \ \ \ \ \ \lim\limits_{x_n
\rightarrow \pm \infty } u ({\rm \mbox{\boldmath $x$}} ) = Z_{\pm}
({\rm \mbox{\boldmath $x$}} ).
\end{equation}
Then at least for $n \leq 8$, the $0$-level set of $u$: $\varGamma:=
\left\{ {\rm \mbox{\boldmath $x$}}\, \left| \,  u ( {\rm
\mbox{\boldmath $x$}} )= 0 \right. \right\} $ has the following
properties:

{\rm (i)} $\varGamma$ is the graph of a function $x_n = h ({\rm
\mbox{\boldmath $x$}}')$ and $h ({\rm \mbox{\boldmath $x$}}')$ is
quasi-periodic if $A,B,Z$ are periodic;

{\rm (ii)} $\varGamma$ is the graph of a function $x_n = h ({\rm
\mbox{\boldmath $x$}}')$ and $h ({\rm \mbox{\boldmath $x$}}')$ is
almost periodic if $A,B,Z$ are almost periodic and if they satisfy
similar conditions as in \eqref{ap-int-1}. }

\vskip 6pt

In the periodic case, denote by $X$ the period of $A,B$ and $Z$ in
$x_n$-direction. One see that ${\frac {\partial u} {\partial x_n}}
>0$ as in De Giorgi's conjecture may be not always true, instead,
$$
u ({\rm \mbox{\boldmath $x$}}', x_n) < u ({\rm \mbox{\boldmath
$x$}}', x_n + X),\ \ \ \ \ \ \ \ \ \ \ \ \ \ \ \ {\rm
\mbox{\boldmath $x$}} \in \mathbb{R}^n
$$
may play a similar role and be used as an additional condition for
the conjecture.

\section{Systems of Reaction-diffusion Equations}

It is well known that the solution $u$ of FitzHugh-Nagumo equation
\begin{equation}\label{FHN-homo}
\left\{
  \begin{array}{ll}
  u_t = \Delta u + \frac{1}{\ve^2} (u-u^3 - v), &\qquad
  {\rm \mbox{\boldmath $x$}}  \in
  \mathbb{R}^n,\ \ t>0, \\
  v_t = \Delta v + u- v, & \qquad
  {\rm \mbox{\boldmath $x$}}  \in
  \mathbb{R}^n,\ \ t>0
\end{array}
\right.
\end{equation}
behaviors as that of \eqref{RD-homo}. From \eqref{FHN-homo} one also
derives the law of the propagation of the level set of $u$, which is
similar to \eqref{MCF-homo} (\cite{Lou, nmhs}). For such reasons, we
propose an analogue of De Giorgi's conjecture for the following
spatially inhomogeneous FitzHugh-Nagumo equations:
\begin{equation}\label{FHN-inhomo}
\left\{
  \begin{array}{ll}
  u_t = \nabla( A({\rm \mbox{\boldmath $x$}}) \nabla u)
   + \frac{1}{\ve^2} B({\rm \mbox{\boldmath $x$}}) u(Z^2 ({\rm \mbox{\boldmath $x$}})
   -u^2 ) - D({\rm \mbox{\boldmath $x$}}) v, &\qquad   {\rm \mbox{\boldmath $x$}}  \in
  \mathbb{R}^n,\ \ t>0, \\
  v_t = \nabla(\alpha({\rm \mbox{\boldmath $x$}}) \nabla v) +
  \beta({\rm \mbox{\boldmath $x$}}) u - \gamma({\rm \mbox{\boldmath $x$}}) v, & \qquad
  {\rm \mbox{\boldmath $x$}}  \in
  \mathbb{R}^n,\ \ t>0.
\end{array}
\right.
\end{equation}

\vskip 6pt

\underline{\bf De Giorgi Type Conjecture for FitzHugh-Nagumo
Equations}. {\it Let $A,$ $ B,D,Z, \alpha, \beta, \gamma$ be
positive, almost periodic functions, $0<\ve \ll 1$. Assume that
$(u,v)$ is an entire stationary solution of {\rm
(\ref{FHN-inhomo})}, and that
\begin{equation}\label{2limits} Z_- ({\rm \mbox{\boldmath $x$}}) \leq
u ({\rm \mbox{\boldmath $x$}}) \leq Z_+ ({\rm \mbox{\boldmath
$x$}}),\ \ \ \ \ \ \lim\limits_{x_n \rightarrow \pm \infty } u ({\rm
\mbox{\boldmath $x$}} ) = Z_{\pm} ({\rm \mbox{\boldmath $x$}} ).
\end{equation}
Then at least for $n \leq 8$, the $0$-level set of $u$:
$\varGamma:=\left\{ {\rm \mbox{\boldmath $x$}}\, \left| \,  u ( {\rm
\mbox{\boldmath $x$}} )= 0 \right. \right\} $ must be the graph of a
function $x_n = h ({\rm \mbox{\boldmath $x$}}')$.

Moreover, function $h ({\rm \mbox{\boldmath $x$}}')$ is a constant
(resp. quasi-periodic, almost periodic) when $A,B,D,Z, \alpha,
\beta, \gamma$ are constant (resp. periodic, almost periodic and
satisfy similar conditions as in \eqref{ap-int-1}).}

\vskip 6pt

{\it Remark}. Both Bernstein's problem and De Giorgi's conjecture
listed above can be extended to more generale cases. For example, $a
= a( {\rm \mbox{\boldmath $x$}}, {\rm \mbox{\boldmath $n$}})$,
$B({\rm \mbox{\boldmath $x$}}) = B({\rm \mbox{\boldmath $x$}}, u
,\nabla u)$, etc.. One can consider other systems, for example, a
three-component reaction-diffusion system arising in gas discharge
phenomena (cf. \cite{NTYU}). Finally, one can consider spatially
homogeneous and/or inhomogeneous analogue of the above Bernstein's
problem and De Giorgi's conjecture in various manifolds.

\vskip 8mm

{\sc\small Department of Mathematics, Tongji University, Shanghai
200092, China}

{\it\small E-mail address}: blou@tongji.edu.cn


\begin{thebibliography}{123456}

\bibitem{Chen}{\sc X.~Chen}, Generation and propagation of interfaces for
    reaction-diffusion equations, {\it J. Differential Equations},
    {\bf 96} (1992), 116-141.

\bibitem{DeG}{\sc E.~De Giorgi}, Convergence problems for functionals
    and operators, in \lq\lq Proc. Internat. Meeting on Recent Methods in
     Nonlinear Analysis" (Rome, 1978) (E. De Giorgi et al., eds.),
     Pitagora, Bologna, 1979, 131-188.

\bibitem{GG}{\sc N.~Ghoussoub} and {\sc C.~Gui}, On De Giorgi's conjecture in
   dimensions $4$ and $5$, {\it Ann. of Math.}, {\bf 157} (2003), 313-334.

\bibitem{HKMN}{\sc D.~Hilhorst, G.~Karali, H.~Matano, K.~Nakashima},
Singular limit of a spatially inhomogeneous Lotka-Volterra
competition-diffusion system, {\it Comm. Partial Differential
Equations}, {\bf 32} (2007), 879--933.

\bibitem{Lou}{\sc B.D.~Lou}, Singular limits of spatially
    inhomogeneous convection-reaction-diffusion equations,{\it J. Stat.
    Phys.}, {\bf 129} (2007), 509-516.


\bibitem{nmhs}{\sc K.I.~Nakamura, H.~Matano, D.~Hilhorst} and {\sc R.~Schatzle},
   Singular limit of a reaction-diffusion equations with
   a spatially inhomogeneous reaction term, {\it J. Stat. Phys.},
    {\bf 95} (1999), 1165-1185.

\bibitem{NTYU}{\sc Y.~Nishiura, T.~Teramoto, X.~Yuan and K.I.~
Ueda}, Dynamics of traveling pulses in heterogeneous media, {\it
Chaos}, {\bf 17}, 037104 (2007).


\bibitem{Sa}{\sc O.~Savin}, Phase transitions: regularity of flat level
    sets, {\it Ann. of Math.}, (to appear).

\end{thebibliography}
\end{document}